\documentstyle[12pt]{article}

\input amssym.def
\input amssym.tex

\def\vv{{\underline{v}}}
\def\tt{{\underline{t}}}
\def\a{\underline{a}}

\def\Z{\Bbb Z}

\def\C{\Bbb C}

\newtheorem{theorem}{Theorem}

\newenvironment{proposition}
{\smallskip\noindent{\bf Proposition\/}.}{\smallskip\par}

\newenvironment{remark}
{\smallskip\noindent{\bf Remark\/}.}{\smallskip\par}
\newenvironment{remarks}
{\smallskip\noindent{\bf Remarks\/}.}{\smallskip\par}
\newenvironment{corollary}
{\smallskip\noindent{\bf Corollary\/}.}{\smallskip\par}

\title{The Alexander polynomial of a plane curve singularity
and the ring of functions on it}
\author{A.Campillo
\and F.Delgado \thanks{First two authors were partially supported by
DGICYT PB97-0471 and by Junta de Castilla y Le\'on:
VA51/97. Address:
University of Valladolid, Dept. of Algebra, Geometry and Topology,
47005 Valladolid, Spain.
E-mail: campillo\symbol{'100}cpd.uva.es, fdelgado\symbol{'100}agt.uva.es}
\and S.M.Gusein--Zade \thanks{Partially supported by grants Iberdrola,
RFBR--98--01--00612 and INTAS--97--1644. Address: Moscow State University,
Dept. of Mathematics and Mechanics, Moscow, 119899, Russia.
E-mail: sabir\symbol{'100}mccme.ru}}
\date{}

\begin{document}
\sloppy

\maketitle

\begin{abstract}
We give two formulae which express the Alexander polynomial
$\Delta^C$
of several variables of a plane curve singularity $C$ in terms of the ring
${\cal O}_{C}$ of germs of analytic functions on the curve. One of them
expresses $\Delta^C$ in terms of dimensions of some factorspaces
corresponding to a (multi-indexed) filtration on the ring ${\cal O}_{C}$.
The other one gives the coefficients of the Alexander polynomial
$\Delta^C$ as Euler characteristics of some explicitly described
spaces (complements to arrangements of projective hyperplanes).

A version of this text has been published in Russian
Mathematical Surveys, v.54 (1999), N 3 (327), p.157--158.
 \end{abstract}

The ring ${\cal O}_{X}$ of germs of holomorphic functions on a germ $X$ of
an analytic set determines $X$ itself (up to analytic equivalence).
Thus all invariants of $X$, in particular, topological ones, can ``be read"
from ${\cal O}_{X}$. There arises a general problem to find expressions for
invariants of $X$ in terms of the ring ${\cal O}_{X}$. The
Alexander polynomial $\Delta^C$ of several variables is a
complete topological invariant of a plane curve singularity
$C\subset (\C^2, $0) (\cite{Y}). A formula
of D.Eisenbud and W.Neumann ({\cite{EN}}) expresses the Alexander
polynomial in terms of an embedded resolution of the curve $C$.
In this note we give two formulae for the Alexander polynomial
directly in terms of the ring of germs of analytic functions on
the curve $C$. One of them expresses the Alexander polynomial
$\Delta^C$ in terms of dimensions of some factorspaces corresponding
to a (multi-indexed) filtration on the ring ${\cal O}_{C}$. The
other one gives the coefficients of the Alexander polynomial
$\Delta^C$ as Euler characteristics of some explicitly described
spaces (complements to arrangements of projective hyperplanes).
It seems to be the first result which describes the coefficients
of the Alexander polynomial (and thus of the zeta--function of the
monodromy) as Euler characteristics of some spaces. Another formula
which expresses the Lefschetz numbers of iterates of the monodromy
(and therefore the zeta--function of it) for a hypersurface
singularity of any dimension in terms of Euler characteristics of
some subspaces of the space of (truncated) arcs is given in a paper
of J.Denef and F.Loeser (xxx-Preprint series, math.AG/0001105).

Let $C$ be a germ of a reduced plane curve at the origin in $\C^2$ and
let $C=\bigcup\limits_{i=1}^{r}C_i$ be its representation as the union of
irreducible components (with a fixed numbering). Let ${\cal O}_{\C^2, 0}$
be the ring of germs of holomorphic functions at the origin in $\C^2$ and let
$\{f=0\}$ ($f\in {\cal O}_{\C^2, 0}$) be an equation of the curve $C$. Let
${\cal O}_C$ be the ring of germs of analytic functions on $C$
($\cong {\cal O}_{\C^2, 0}/(f)$), and let
$\Delta^C(t_1, \ldots, t_r)$ be the Alexander polynomial of the link
$C\cap S_\varepsilon^3\subset S_\varepsilon^3$ for $\varepsilon > 0$
small enough (see, e.g., {\cite{EN}}).

\begin{remarks}
1. According to the definition, the Alexander polynomial
$\Delta^C(t_1, \ldots, t_r)$
is well defined only up to multiplication by monomials $\pm\tt^{\underline{m}}=
\pm t_1^{m_1}\cdot\ldots\cdot t_r^{m_r}$ ($\tt=(t_1, \ldots, t_r)$,
${\underline{m}}=(m_1, \ldots, m_r)\in \Z^r$). We fix the Alexander
polynomial assuming that it is really a polynomial (i.e., it does not
contain variables with negative powers) and $\Delta^C(0, \ldots, 0)=1$.

\smallskip\noindent 2. There is some difference in definitions (or rather
in descriptions) of the Alexander polynomial for a curve with one branch ($r=1$)
or with many branches ($r>1$) (see, e.g., {\cite{EN}}). In order to have all
the results (Theorems {\ref{theo1}} and {\ref{theo2}} below) valid for $r=1$
as well, for an irreducible curve $C$, $\Delta^C(t)$ should be not
the Alexander polynomial, but rather the zeta-function $\zeta_C(t)$
of the monodromy, equal to the Alexander polynomial
divided by $(1-t)$. In this case $\Delta^C(t)$ is not a polynomial, but an
infinite power series. However for uniformity of the statements we shall
use the name "Alexander polynomial" for this $\Delta^C(t)$ as well.
\end{remarks}

Let $\varphi_i:(\C_i, 0)\to(\C^2, 0)$ be parametrizations (uniformizations)
of the components $C_i$ of the curve $C$, i.e., germs of analytic
maps such that
${\rm{Im}}\,\varphi_i=C_i$ and $\varphi_i$ is an isomorphism between $\C_i$
and $C_i$ outside of the origin.
For a germ $g\in{\cal O}_{\C^2, 0}$, let $v_i=v_i(g)$ and $a_i=a_i(g)$ be the power
of the leading term and the coefficient at it in the power series decomposition
of the germ $g\circ\varphi_i:(\C_i,0)\to \C$ :
$g\circ\varphi_i(t_i)=a_i\cdot t_i^{v_i}+{~terms~of~higher~degree~}$ ($a_i\ne 0$).
If $g\circ\varphi_i(t)\equiv 0$, $v_i(g)$ is assumed to be equal to $\infty$
and $a_i(g)$ is not defined. The numbers $v_i(g)$ and $a_i(g)$ are defined for
elements $g$ of the ring ${\cal O}_C$ of functions on the curve
$C$ as well.

The semigroup $S=S_C$ of the plane curve singularity $C$ is the subsemigroup of
$\Z_{{\ge 0}}^r$ which consists of elements of the form
$\vv(g)=(v_1(g), \ldots, v_r(g))$ for all germs $g\in {\cal O}_C$ with
$v_i(g)<\infty$; $i=1, \ldots, r$. The extended semigroup $\hat S=\hat S_C$
of the plane curve singularity $C$ is the subsemigroup of
$\Z_{{\ge0}}^r\times(\C^*)^r$ which consists of elements of the form
$(\vv(g);\a(g))=(v_1(g), \ldots, v_r(g);a_1(g), \ldots, a_r(g))$ for all
germs $g\in {\cal O}_C$ with $v_i(g)<\infty$, $i=1, \ldots, r$ ({\cite{CDG1}}).

It is known that both the semigroup $S_C$ and the Alexander polynomial
$\Delta^C(t_1, \ldots, t_r)$ are complete topological invariants of a
plane curve singularity, i.e., each of them determines the germ $C$ up
to topological equivalence ({\cite{W}}, {\cite{Y}}).
Therefore it is interesting to understand a connection between them.
In fact from the formula for the Alexander polynomial in terms of a
resolution of a plane curve singularity (see {\cite{EN}}) it is not difficult
to understand that the Alexander polynomial $\Delta^C(t_1, \ldots, t_r)$
may contain with non-zero coefficients only monomials $\tt^\vv$ for $\vv$
from the semigroup $S_C$ of the curve $C$. For the case of an irreducible
curve $C$ ($r=1$) the corresponding connection has been described in
{\cite{CDG2}}. In this case
$$
\zeta_C(t)=\sum\limits_{i\in S_C}t^i
$$
($S_C\subset\Z_{\ge 0}$).

Let $\pi:\hat S_C\to \Z^r$ be the natural projection: $(\vv, \a)\mapsto \vv$.
For an element $\vv\in \Z^r$, let
$F_{\vv}=\pi^{-1}(\vv)\subset\{\vv\}\times(\C^*)^r\subset\{\vv\}\times\C^r$
be the corresponding fibre of the extended semigroup ({\cite{CDG1}}).
The fibre $F_{\vv}$ is not empty if and only if $\vv\in S_C$.
For $\vv=(v_1, \ldots, v_r)\in \Z^r$, let $J(\vv)=\{g\in {\cal O}_C:
v_i(g)\ge v_i; i=1, \ldots, r\}$ be an ideal in ${\cal O}_C$. One has a natural
linear map $j_{\vv}:J(\vv)\to\C^r$, which sends $g\in J(\vv)$ to
$(a_1, \ldots, a_r)$, where $a_i$ is the coefficient in the power series
expansion $g\circ\varphi_i(t_i)=a_it_i^{v_i}+\ldots$ (the number $a_i$
may be equal to zero). Let $C(\vv)\subset \C^r$ be the image of the map
$j_{\vv}$, let $c(\vv)=\mbox{dim }C(\vv)$. It is not difficult to see that
$C(\vv)\cong J(\vv)/J(\vv+{\underline{1}})$, where ${\underline{1}}=(1, \ldots, 1)$,
and that $F_{\vv}=C(\vv)\cap(\C^*)^r$ (under the natural identification of
$\{\vv\}\times(\C^*)^r$ and $(\C^*)^r$). Therefore the fibre $F_{\vv}$
($\vv\in S_C$) is the complement to an arrangement of linear
hyperplanes in the vector space $C(\vv)$. The extended semigroup
$\hat S_C$ contains some analytic information about the plane curve
singularity $C$, however the dimensions $c(v)$ depend only on the
topological type of $C$ (see {\cite{CDG1}}).

Let ${\cal L}=\Z[[t_1, \ldots, t_r, t_1^{-1}, \ldots, t_r^{-1}]]$ be
the set of formal Laurent series in $t_1, \ldots, t_r$. Elements of ${\cal L}$
are expressions of the form $\sum\limits_{\vv\in\Z^r}k(\vv)\cdot\tt^\vv$ with
$k(\vv)\in\Z$, generally speaking, infinite in all directions. ${\cal L}$
is not a ring, but a $\Z[t_1, \ldots, t_r]$-- (or even
$\Z[t_1, \ldots, t_r, t_1^{-1}, \ldots, t_r^{-1}]$--) module. The polynomial
ring $\Z[t_1, \ldots, t_r]$ can be in a natural way considered as being embedded
into ${\cal L}$.

Let $L_C(t_1, \ldots, t_r)=\sum\limits_{\vv\in\Z^r}c(\vv)\cdot\tt^\vv\ \in\ {\cal L}$,
$P'_C(t_1, \ldots, t_r)=(t_1-1)\cdot\ldots\cdot(t_r-1)\cdot L_C(t_1, \ldots, t_r)$.
One can easily see that $P'_C(t_1, \ldots, t_r)$ is in fact a polynomial, i.e.,
$P'_C(t_1, \ldots, t_r)\in \Z[t_1, \ldots, t_r]$. This follows
from the fact that, if $v_i'$ and $v_i''$ are negative, then
$c(v_1, \ldots, v_i', \ldots, v_r)=c(v_1, \ldots, v_i'', \ldots, v_r)$.
Let $P_C(t_1, \ldots, t_r)=P'_C(t_1, \ldots, t_r)/(t_1\cdot\ldots\cdot t_r-1)
\in \Z[[t_1, \ldots, t_r]]$.

\begin{proposition}
For $r>1$, the polynomial $P'_C(t_1, \ldots, t_r)$ is divisible by
$(t_1\cdot\ldots\cdot t_r-1)$, i.e., $P_C(t_1, \ldots, t_r)
\in \Z[t_1, \ldots, t_r]$.
\end{proposition}

For $r=1$, $P_C(t)=L_C(t)$.

\begin{theorem}\label{theo1}
$P_C(t_1, \ldots, t_r)=\Delta^C(t_1, \ldots, t_r)$.
\end{theorem}

The fibre $F_{\vv}$ of the extended semigroup is invariant with respect to
multiplication by non-zero complex numbers. Let ${\Bbb P}(F_{\vv})$ be the
projectivization of the fibre $F_{\vv}$, i.e., ${\Bbb P}(F_{\vv})=F_{\vv}/\C^*$.
The projectivization ${\Bbb P}(F_{\vv})$ of the fibre $F_{\vv}$
is the complement to an arrangement of projective hyperplanes in
a projective space. If
$\vv\ge {\underline{\delta}}$, where ${\underline{\delta}}$ is the conductor of
the semigroup $S_C$ of the curve $C$, then the fibre $F_{\vv}$ coincides with
$(\C^*)^r$ and the Euler characteristic $\chi({\Bbb P}(F_{\vv}))$ of its
projectivization is equal to $1$ for $r=1$ and to $0$ for $r>1$.
Let $\chi({\Bbb P}\hat S_C):=
\sum\limits_{\vv\in\Z_{{\ge0}}^r}\chi({\Bbb P}(F_{\vv}))\cdot\tt^\vv$.

\begin{theorem}\label{theo2}
$$\Delta^C(t_1, \ldots, t_r)=\chi({\Bbb P}\hat S_C).\eqno(*)$$
\end{theorem}

Let $\zeta_C(t)$ ($=\Delta^C(t, t, \ldots, t)$) be the
zeta--function of the monodromy of the germ $f$ (the equation of
the curve $C$). Let $\vert\vv\vert:=v_1+\ldots+v_r$.

\begin{corollary}
$\zeta_C(t)=
\sum\limits_{i=0}^\infty\chi\left(\bigcup\limits_{\vv:\vert\vv\vert=i}{\Bbb
P}(F_{\vv})\right)\cdot\tt^\vv$.
\end{corollary}

\begin{remark}
For an irreducible plane curve singularity all coefficients of
the zeta--function of the monodromy are equal to $0$ or $1$. In
terms of the equation ($*$), $0=\chi(\emptyset)$,
$1=\chi(point)$.
\end{remark}

The proof consists of calculation of the polynomial
$\chi({\Bbb P}\hat S_C)$ in terms of a suitable (not minimal
one)embedded resolution of the curve $C\subset(\C^2, 0)$
and comparing it with the formula for the Alexander polynomial
from {\cite{EN}}. These calculations involve a detailed
knowledge about the structure of the semigroup and its relation
with the resolution of a singularity. In fact the polynomials
$P_C(t_1, \ldots, t_r)$ and $\chi({\Bbb P}\hat S_C)$ coincide
for any (not necesseraly plane) curve.
The proof will be published elsewhere.

A global version of the result from {\cite{CDG2}} for a plane
algebraic curve with one place at infinity was obtained in
{\cite{CDG3}}.

\end{document}